\documentclass[a4paper,10pt]{amsart}
\usepackage[utf8]{inputenc}
\usepackage{amsmath}

\begin{document}

\title{On the variety of four dimensional Lie algebras}
\author{Laurent Manivel}
\date{February 2015}

\newtheorem{theorem}{Theorem}[section]
\newtheorem{lemma}[theorem]{Lemma}
\newtheorem{proposition}[theorem]{Proposition}
\newtheorem{corollary}[theorem]{Corollary}
\newtheorem{definition}[theorem]{Definition}

\def\proof{\noindent {\it Proof}.\hspace{1mm} }
\def\PP{\mathbf{P}}
\def\cA{{\mathcal A}}
\def\cB{{\mathcal B}}
\def\cO{{\mathcal O}}
\def\cL{{\mathcal L}}
\def\ra{\rightarrow}\def\lra{\longrightarrow}
\def\fg{\mathbf{g}}\def\fh{\mathbf{h}}\def\fk{\mathbf{k}}
\def\fsl{{\mathbf sl}}\def\fgl{{\mathbf gl}}
\def\s{{\sigma}}\def\t{{\tau}}\def\iot{{\iota}}
\def\fine1{\mathbf{aff}_1}
\def\affine2{\mathbf{aff}_2}
\def\bz{\mathbf{z}}
\def\he{\mathbf{he}_3}
\def\abel{\mathbf{ab}}
\def\ab1{\mathbf{ab}_1}
\def\ab2{\mathbf{ab}_2}
\def\Lie{\mathrm{Lie}}

\begin{abstract}  Lie algebras of dimension $n$ are defined by their structure constants,  which can be seen
as sets of $N=n^2(n-1)/2$ scalars (if we take into account the skew-symmetry condition) to which the Jacobi
identity imposes certain quadratic conditions. Up to rescaling, we can consider such a set as a point in the projective
space $\PP^{N-1}$. Suppose $n=$4, hence $N=24$. Take a random subspace of dimension $12$ in $\PP^{23}$,
over the complex numbers. 
We prove that this subspace will contain exactly $1033$ points giving the structure constants of some 
four dimensional Lie algebras. Among those, $660$ will be isomorphic to $\fgl_2$, $195$ will be the sum of two 
copies of the Lie algebra of one dimensional affine transformations, $121$ will have an abelian,
three-dimensional  derived algebra, and $57$ will have for derived algebra
the three dimensional Heisenberg algebra. This answers a question of Kirillov and Neretin.    
\end{abstract}

\maketitle

\section{The variety of Lie algebras}

A Lie algebra structure on a $k$-vector space $V_n$ of finite dimension $n$ is given 
by a Lie bracket, which can be considered as a linear map from $\wedge^2V_n$ to 
$V_n$, that we will denote by $\omega$:
$$  \omega(X\wedge Y)=[X,Y] \qquad\forall X,Y\in V_n.$$
A Lie bracket needs to verify the Jacobi identity
$$[[X,Y],Z]+[[Y,Z],X]+[[Z,X],Y]=0 \qquad\forall X,Y,Z\in V_n.$$
The left hand side of this identity is skew-symmetric in the three arguments,
and what is required in terms of $\omega$ is  the vanishing of the map $Jacobi(\omega)$
from $\wedge^3V_n$ obtained as the composition
$$\wedge^3V_n\hookrightarrow \wedge^2V_n\otimes V_n\stackrel{\omega\otimes id}{\lra}
V_n\otimes V_n\lra \wedge^2V_n\stackrel{\omega}{\lra} V_n.$$
The variety of $n$-dimensional Lie algebras is 
$$\Lie_n=\{\omega\in Hom(\wedge^2V_n, V_n), \quad Jacobi(\omega)=0\}.$$
Of course this is a cone, with vertex the point defining the abelian Lie algebra $\abel_n$. 
So we do not loose much information by considering instead the projective variety
$\PP\Lie_n$. This is a projective variety defined by a collection of 
quadratic equations, invariant under the natural action of $GL(V_n)$. So it is 
natural to ask: 
\begin{enumerate}
 \item What are the irreducible components of $\PP\Lie_n$?
 \item What is the dimension, the degree, more generally the geometry of each component? 
 \item What is the Lie algebra parametrized by the generic point of each component?
\end{enumerate}
These questions have been addressed in \cite{kn}, which contains a number of interesting 
results. For example, the authors prove that the number of irreducible components in 
$\PP\Lie_n$ grows at least exponentially in the square root of $n$. Moreover
they give the list of these components for $n\le 6$, with their dimensions and the Lie
algebras corresponding to their generic points. They do not compute their degrees, 
and in fact they conclude their paper by asking what are the degrees of the four 
components of $\PP\Lie_4$. This is the question we want to answer in this note. 

\section{Classification of low dimensional Lie algebras}

From now on we work over the complex numbers. 
The classification of Lie algebras is well known in dimension four \cite{bs, lauret}
(for more information in higher dimension, see for example \cite{yam}). 
We will need some details on this classification in order 
to construct smooth models of the components of $\PP\Lie_4$ in the next section. 

\subsection{Classification in dimension smaller than four}

In dimension $n=2$, a Lie bracket is defined as $[X,Y]=\theta(X,Y)V$, where 
$\theta$ is a fixed non zero skewsymmetric form on $V_2$ and the vector $V$ 
can be chosen arbitrarily. This means that $\Lie_2\simeq V_2$, 
where the origin corresponds to the abelian Lie algebra structure. For any 
other vector, the resulting algebra is isomorphic to the Lie algebra of 
the group of affine tranformations of the affine line, that we denote by $\fine1$.
We conclude that $\PP\Lie_2\simeq \PP^1$ with its usual transitive action 
of $PGL_2$. 

\smallskip
In dimension $n=3$, a Lie bracket is defined by $\omega\in Hom(\wedge^2V_3, V_3)$. 
If we fix a non degenerate trilinear non zero skewsymmetric form $\theta$ on $V_3$, it identifies
 $\wedge^2V_3$ with $V_3^*$, and we 
get an induced identification of $Hom(\wedge^2V_3, V_3)$ with 
$$V_3\otimes V_3=Sym^2V_3 \oplus V_3^*.$$
Therefore $\omega$ can be interpreted as $q+u$, with $q$ a symmetric tensor and $u$
a linear form. Moreover, a computation shows that the Jacobi identity translates 
into the condition that the contraction of $q$ with $u$ is zero. One can then 
distinguish two cases. Either $u=0$ and $q$ can be arbitrary, or $u$ is non zero
and $q$ must be a degenerate tensor, with $u$ contained in its kernel. This 
gives two irreducible components in $\PP\Lie_3$. The first one is the 
linear space $C_1=\PP (Sym^2V_3)\simeq\PP^5$. The second one is 
$$C_2=\PP \{q+Ker(q), \quad q\in Sym^2V_3, \quad \det(q)=0\}.$$
It is birationally a one-dimensional bundle over the determinantal cubic fourfold 
parametrizing singular tensors in  $\PP(Sym^2V_3)$, which is also the intersection 
of the two components. 

In order to be more intrinsic, we will denote by $Hom_s(\wedge^2V_3, V_3)\subset 
Hom(\wedge^2V_3, V_3)$ the space corresponding to symmetric tensors and defining 
$C_1$. Generically, such a tensor has rank three and the associated Lie algebra 
is then isomorphic to $\fsl_2$. When the rank is one the associated Lie algebra
is the Heisenberg algebra $\he$  (denoted $\Gamma_3$   in \cite{kn}). 

\subsection{Classification in dimension four}

The first observation is that if $\fg$ is a four dimensional Lie algebra,
the derived algebra $\fh=[\fg,\fg]$ is a proper ideal of $\fg$. 
This allows to proceed by induction. In fact $\fh$ cannot be any Lie algebra 
of dimension at most three. A routine examination shows that the only possibilities
are $\fh = \fsl_2, \he, \abel_3, \abel_2, \abel_1$ or $0$. Moreover the first four 
possibilities define the four irreducible components of $\Lie_4$ or $\PP\Lie_4$.

\begin{proposition}
The variety $\PP\Lie_4$ has four irreducible components, all of dimension $11$.
They are characterized by the following properties of the Lie algebras defined by their generic points:
\begin{enumerate}
\item  $C_1$: be isomorphic to $\fgl_2$,
\item  $C_2$: have derived algebra isomorphic to $\he$,
\item  $C_3$: have derived algebra isomorphic to  $\abel_3$,
\item  $C_4$: be isomorphic to $2\fine1$.
\end{enumerate}
\end{proposition}

Note that the irreducible components of $Lie_n$ have been classified up to $n=7$
\cite{cd}. Starting from $n=5$ there exist components of different dimensions. 

\section{Smooth models of the irreducible components}

Desingularizations of the components of $\PP\Lie_4$ have been described in \cite{basili}, 
and used to describe their singular loci. In this section we give more intrinsic descriptions 
of smooth models of the components, and use them to compute their degrees. 

\subsection{The first component}

Suppose that $V_4$ is endowed with a Lie algebra structure isomorphic to $\fgl_2$. Let $U\subset V_4$ be the derived
algebra, isomorphic to $\fsl_2$. For any $T\in V_4$, $ad(T)$ restricts to a derivation of $U$, so there is a unique vector
$p(T)\in U$ such that $ad(T)_{|U}=ad(p(T))$. Or course the restriction of $p$ to $U$ must be the identity. As we have seen,
the  $\fsl_2$-structure on $U$ is defined by a non-degenerate symmetric tensor $\sigma\in Hom_s(\wedge^2U,U)$. 

Conversely, let $G_U$ be the subspace of $Hom(V_4,U)$ consisting of homomorphisms whose restriction to $U$ 
is a multiple of the identity. Let $H_U=Hom_s(\wedge^2U,U)$.  A pair $(p,\sigma)$ with $p\in G_U$ and 
$\sigma\in H_U$  defines a Lie algebra structure on $V_4$, the corresponding bracket 
$\omega\in Hom(\wedge^2V_4,V_4)$ being given by 
$$\omega(T,X)=\sigma(p(T),X) \qquad \forall T\in V_4, \; X\in U.$$
Note that $G_U$ and $H_U$ are the fibers of two vector bundles $G$ and $H$ over $\PP(V_4^*)$, of respective 
ranks $4$ and $6$. The previous discussion implies:
 
\begin{proposition}
There exists a  birational map $$\pi_1: \PP (G)\times _{\PP(V_4^*)}\PP (H)\longrightarrow \PP C_1 .$$
\end{proposition}

Beware that $\pi_1$ is not defined on the locus $\Sigma$ where $p$ vanishes on $U$, and the line $p(V_4)$
is contained in the kernel of $\sigma$. This is a fibered version of the usual desingularisation of the symmetric
determinantal cubic in $\PP^5$. In particular $\Sigma$ is smooth, and we would need to blow it up in order
to regularize $\pi_1$. 

The degree of $C_1$, as a subvariety of $\PP^{23}$, can be computed by integrating the first Chern class 
of the hyperplane line bundle $O(1)$, taken to the power $\dim C_1 =11$. This integration can be 
performed on any birational model of $C_1$. In particular, since  $\pi_1^*O(1)=O_{G,H}(1,1)$, even if $\pi_1$ is not a morphism we can compute the degree of $ \PP C_1$ as 
$$deg  \;\PP C_1=\int_{ \PP C_1}c_1(O(1))^{11}=\int_{ \PP (G)\times _{\PP(V_4^*)}\PP (H)}c_1(O_{G,H}(1,1))^{11}.$$
Recall that if $E$ is a vector bundle of rank $e$ on a variety $X$, the projective bundle $\PP(E)$ is endowed 
with a tautological line bundle $O_E(-1)$. If $p : \PP(E)\rightarrow X$ denotes the projection, $O_E(-1)$ is a sub-bundle 
of $p^*E$. Let $O_E(1)$ be its dual. Then the Segre classes of $E$ can be obtained by push-forward:
$$p_* c_1(O_E(1))^k=s_{k-e+1}(E),$$
with the convention that the right hand side is zero when $k\le e-2$.
Alternatively, the total Segre class $s(E)=\sum_{j\ge 0}s_j(E)$ is the formal inverse of the total Chern class 
\cite[Chapter 3]{fulton}. 

Considering  $\PP (G)\times _{\PP(V_4^*)}\PP(H)$ as an iterated projective bundle, and pushing-forward twice,
we obtain 
$$deg  \;\PP C_1=\sum_{i+j=11}\frac{11!}{i!j!}\int_{\PP(V_4^*)}s_{i-3}(G)s_{j-5}(H).$$
From the tautological exact sequence on $\PP(V_4^*)$ we get that 
$$s(G)=(1-h)^{-4}, \qquad s(H)=(1+2h)^4(1+h)^{-10},$$
in terms of the hyperplane class $h$. A  computation yields:

\begin{proposition}
The degree of $\PP C_1$ is $660$.
\end{proposition}
 
\smallskip 

\subsection{The second component}

Suppose that $V_4$ is endowed with a Lie algebra structure such that the derived algebra is isomorphic to the three
dimensional Heisenberg algebra $\he$. Denote by $U$ this hyperplane of $V_4$. The Heisenberg algebra structure on
$U$ is determined, up to scalars, by the one-dimensional center $L\subset U$, which is also the derived algebra of $U$. 

Now let $T$ be an arbitrary element of $V_4$. Then $ad(T)$ restricts to an endomorphism $\theta(T)$ of $U$ which 
preserves $L$. We thus get two induced endomorphisms  $\theta_{U/L}(T)$ of $U/L$ and  $\theta_{L}(T)$ of $L$. Of 
course the latter is just a scalar. 

\begin{lemma}
The Jacobi identity is equivalent to the condition that 
$$ \theta_{L}(T)=\mathrm{trace}\; \theta_{U/L}(T) \qquad\forall T\in V_4. $$
\end{lemma}

Conversely, let  $\Omega$ belong to $Hom_L(\wedge^2U,L)$, the one dimensional space of two-forms whose kernel contains $L$. Denote by $End_L(U)$ the vector space of endomorphisms of 
$U$ preserving $L$, and by  $End^0_L(U)$ the hyperplane defined by the condition 
of the Lemma. Let $\theta\in Hom (V_4, End^0_L(U))$. Choose some $T\notin U$. Then for any
 $\alpha\in\wedge^2V_4$, we can find vectors $A,B,C$ in $U$ such that $\alpha=T\wedge A+B\wedge C$. 
We would therefore like to define our Lie algebra structure on $V_4$ by letting 
$$\omega(\alpha)=\theta(T)(A)+\Omega(B,C).$$
For this to make sense, it must be independent of $T$. This is clearly equivalent to the condition that $\theta(D)(A)=
\Omega(D,A)$ when $D,A$ belong to $U$. 

Consider therefore the vector space $F_{L,U}$ of pairs $(\Omega,\theta)$, where $\Omega$ belongs to  $Hom_L(\wedge^2U,L)$, 
$\theta$ to $End^0_L(U)$, and the previous compatibility condition is satisfied. This defines a rank seven vector bundle on the flag variety $F(1,3,V_4)$, which has dimension $5$. 

\begin{proposition}
There exists a proper birational map $\pi_2: \PP (F)\lra \PP C_2$.
\end{proposition}

\proof The map $\pi_2$ has just been described. It is clearly birational since generically, $U$ can only be the derived algebra, 
$L$ its center, and the pair $(\Omega,\theta)$ is completely determined by the Lie bracket. $\Box$

\medskip
Since obviously $\pi_2^*O(1)=O_F(1)$, we can obtain the degree of $ \PP C_2$ as 
$$deg \; \PP C_2=\int_{ \PP C_2}c_1(O(1))^{11}=\int_{ \PP (F)}c_1(O_F(1))^{11}=\int_{F(1,3,V_4)}s_5(F).$$

In order to compute this degree, first observe that 
$$End^0_L(U)\simeq Hom(U/L,U)\qquad \mathrm{and} \qquad Hom_L(\wedge^2U,L)\simeq Hom(\wedge^2(U/L),L).$$
Denote the latter line bundle by $M$. Mapping the pair $(\Omega,\theta)$ to $\Omega$ defines a vector 
bundle morphism from $F$ to $M$ with kernel $K\simeq Hom(V_4/U, End^0_L(U)).$
From the exact sequence
$$0\lra  End^0_L(U)\lra End(U)\lra Hom(L,U)\lra 0$$
we deduce that the total Segre class of $F$ is 
$$s(F)=s(M)s( Hom(V_4/U, End(U)))c( Hom(V_4/U, Hom(L,U))).$$
Let us denote by $p$ the projection from $F(1,3,V_4)$ to $\PP (V_4^*)$, which identifies the flag variety with the
projective bundle $\PP (U)$. The tautological bundle $O_U(-1)$ is our $L$. 
 The hyperplane line bundle on $\PP (V_4^*)$ is $O(-1)=V_4/U$. Using the projection formula, we get that 
$$p_*s(F)=s(End(U)(-1))p_*S, $$
where $S=s(M)c(U(-1))\otimes L^*)$. 

\begin{lemma}
The push-forward of $S$ is
$$p_* S=13-90h+318h^2-738h^3.$$
\end{lemma}

\proof In order to compute this push-forward, we need to decompose $S$ as
$S=\sum_k \ell^kp^*S_k$, where $\ell=c_1(L^*)$, and use the fact that $p_*\ell^k=s_{k-2}(U)$. Note that from the 
tautological sequence on  $\PP (V_4^*)$, the total Segre class of $U$ is $s(U)=1+h$. Hence $p_*S=S_2+hS_3$.
Let $v_1,v_2,v_3$ be the Chern roots of $U(-1)$. Then  
 $$c(U(-1))\otimes L^*)=\prod_{i=1}^3(1+v_i+\ell)=\sum_{j=0}^3e_j\ell^{3-j}.$$
Here we denoted by $e_j$ the degree $j$ elementary symmetric function of   $1+v_1,1+v_2,1+v_3$. 
Since $c(M)=1-2\ell+h$, we get that
$$S_k=\sum_{p,j}\frac{(p+k+j-3)!}{p!(k+j-3)!}(-1)^p 2^{k+j-3}e_jh^p.$$
The sum over $p$ gives $(1+h)^{-(k+j-2)}$, and we readily deduce that
$$p_*S=\sum_{j\ge 1}\frac{2^{j-1}e_{j}}{(1+h)^{j}}+h\sum_{j\ge 0}\frac{2^je_j}{(1+h)^{j+1}}.$$
 In order to evaluate that sum, it is convenient to introduce the stretched variant of the total Chern class
$$c_z(U(-1)):=\sum_{j\ge 0}z^jc_j(U(-1))=(1-zh)^{4}.$$
where the last equality follows from the twisted tautological exact sequence. A formal substitution  yields 
$$e_z(U(-1)):=\sum_{j\ge 0}z^je_j=(1+z)^3c_{\frac{z}{1+z}}(U(-1))=(1+z-zh)^{4}/(1+z).$$
We can then get $p_*S$ as 
$$p_*S=\frac{1}{2}\Big(e_{\frac{2}{1+h}}(U(-1))-1\Big)-\frac{h}{1+h}e_{\frac{2}{1+h}}(U(-1)),$$
whence the result. $\Box$

\medskip
In order to deduce the push-forward of $s(F)$, there remains to observe that, again from the tautological exact 
sequence, we can deduce that 
$$c(End(U)(-1))=\frac{1-h}{(1+h)^4}.$$
The final result is the following:

\begin{proposition}
The degree of $\PP C_2$ is $57$.
\end{proposition}

\medskip 

\subsection{The third component}

Suppose $V_4$ is endowed with a Lie algebra structure such that the derived algebra is abelian of dimension three. 
More generally we could ask that there exists a hyperplane $U$ of $V_4$ over which the Lie bracket vanishes 
identically. In order to define the bracket completely, what remains to do is to prescribe the bracket of an element
of $U$ with a given element not in $U$. More intrinsically, we need to specify a map
$$  \theta : V_4/U \lra End(U),$$
which can be arbitrary: the Jacobi identity will always hold true. 
As before, denote by $U$ the tautological rank three vector bundle on $\PP(V_4^*)$. 
Let $E$ be the rank $9$ vector bundle $Hom (V_4/U \, End(U))$ on  $\PP(V_4^*)$.

\begin{proposition}
There exists a proper birational map $\pi_3: \PP (E)\lra \PP C_3$.
\end{proposition}

\proof What remains to be proved is that the natural map $\pi_3: \PP (E)\lra \PP C_3$ we have just described is 
birational. This follows from the obvious fact that the general point of $C_3$ defines a Lie algebra whose unique 
abelian codimension one subalgebra is the derived subalgebra. $\Box$

\medskip 
In particular we can obtain the degree of $ \PP C_3$ as 
$$deg  \;\PP C_3=\int_{ \PP C_3}c_1(O(1))^{11}=\int_{ \PP (E)}c_1(O_E(1))^{11}=\int_{\PP(V_4^*)}s_3(E).$$
In order to compute this number it is convenient to use the Chern character. Indeed we have $E=U^*\otimes 
\wedge^2U^*$ and the tautological exact sequence yields, if we denote by $h$ the hyperplane class on 
$\PP(V_4^*)$,
$$ch(U^*)=4-e^{-h}, \qquad ch(\wedge^2U^*)=6-4e^{-h}+e^{-2h}.$$
The Chern character being multiplicative, $ch(E)=ch(U^*) ch(\wedge^2U^*)$, from which we compute that 
$p_1(E)=9h$, $p_2(E)=h^2$, $p_3(E)=-15h^3$, and finally 
$$s_3(E)=\frac{1}{6}(p_1(E)^3+3p_1(E)p_2(E)+2p_3(E))=121h^3.$$
We have proved that:

\begin{proposition}
The degree of $\PP C_3$ is $121$.
\end{proposition}
 
 \smallskip
 
\subsection{The fourth component}

Suppose $V_4$ is endowed with a Lie algebra structure such that the derived algebra is abelian and two-dimensional. 
Denote the derived algebra by $U$. The Lie bracket $\omega\in Hom(\wedge^2V_4,V_4)$ factorizes through $U$ and 
vanishes on $\wedge^2U$. Since $V_4\wedge U / \wedge^2U\simeq (V_4/U) \otimes U$, the restriction of $\omega$ 
to  $V_4\wedge U$ induces a morphism 
$$ \tau : V_4/U \lra End(U).$$
This morphism does not define completely the Lie algebra structure, but observe that $\wedge^2V_4 / V_4\wedge U
\simeq \wedge^2(V_4/U)$ is one dimensional, so the data that is missing is just a skew-symmetric $U$-valued two-form
on $V_4/U$. Moreover, the Jacobi identity will not involve this part of the Lie bracket, but only $\tau$:

\begin{lemma}
The Jacobi identity is equivalent to the vanishing of the composition 
$$ \wedge^2(V_4/U) \stackrel{\wedge^2\tau}{\lra} \wedge^2 End(U)
\stackrel{ [\, ,\; ]}{\lra} \fsl (U).$$
\end{lemma}

Obviously this only depend on the traceless part $\tau_0$ of $\tau$, and the Jacobi identity imposes that $\tau_0$
belongs to the Segre product
$$\PP (V_4/U)^*\times \PP(\fsl(U))\subset\PP (Hom(V_4/U, \fsl(U))).$$


So let $S$ denote this relative Segre product over the Grassmannian $G(2,V_4)$: this is a smooth fiber bundle 
with fiber $\PP^1\times \PP^2$. Let $s$ denote the projection to $G(2,V_4)$. We define a vector bundle $P$ 
over $S$ as follows: the fiber of $P$ above the line 
$[\tau_0]$ is the vector space of  Lie brackets $\omega\in Hom(\wedge^2V_4,U)\subset Hom(\wedge^2V_4,V_4)$
that vanish on $\wedge^2U$, and such that the traceless part of the associated morphism $\tau$ is a multiple 
of $\tau_0$. By mapping $\omega$ to $\tau$ we obtain the exact sequence of vector bundles 
$$0\lra s^*Hom(\wedge^2(V_4/U), U)\lra P\lra O_S(-1)\oplus s^*(V_4/U)^*\lra 0.  $$

\begin{proposition}
There exists a proper birational map $\pi_4: \PP (P)\lra \PP C_4$.
\end{proposition}

\proof The map $\pi_4$ has just been described. It is clearly birational since generically, $U$ can only be the derived algebra, 
and the Lie bracket is exactly a point in the fiber of $P$ over $U$. $\Box$

\medskip
Here again $\pi_4^*O(1)=O_P(1)$, and  we can obtain the degree of $ \PP C_4$ as 
$$deg \; \PP C_4=\int_{ \PP C_4}c_1(O(1))^{11}=\int_{ \PP (P)}c_1(O_P(1))^{11}=\int_S s_7(P).$$
Observe that in order to compute this number, we may suppose that the exact sequence above is split, 
hence that $P\simeq  O_S(-1)\oplus s^*R$ with $R=Q^*\oplus U(-1)$. Recall that, like the total Chern class,
the total Segre class is multiplicative with respect to direct sums, so that, using the projection formula, we get 
$$s_* (s(P))= s_* (s(O_S(-1))) s(R).$$
Since $S$ is a bundle of products of projective spaces over the Grassmannian, the usual formula for 
projective bundles gives 
$$s_*(c_1(O_S(1))^k)=\sum_{i+j=k}\frac{k!}{i!j!}s_{i-1}(Q)s_{j-2}(\fsl(U)),$$
and therefore we get the degree of $\PP C_4$ as
$$deg \; \PP C_4=\sum_{i,j}\frac{(i+j)!}{i!j!}\int_{ G(2,V_4)}s_{i-1}(Q)s_{j-2}(\fsl(U))s_{7-i-j}(R).$$
This computation is a straightforward excercise in elementary Schubert calculus. Note in particular that $\fsl(U)$ being 
self-dual has zero Chern and Segre classes in odd degrees. More precisely, in terms of the Chern classes $c_1$ and
$c_2$ of $U^*$ we get 
$$c(\fsl(U))=1+4c_2-c_1^2, \qquad s(\fsl(U))=1-4c_2+c_1^2+10c_2^2.$$
The Segre classes of the two other bundles are easily computed:
$$s(Q)=1-c_1+c_1^2-c_2, \qquad s(R)=1+4c_1+10c_1^2+40c_1c_2+70c_2^2.$$
Keeping in mind that $c_2^2=c_2c_1^2=1$ is the punctual class, while $c_1^4=2$, 
we finally obtain that:

\begin{proposition}
The degree of  $\PP C_4$ is equal to $195$.
\end{proposition}

\bigskip\noindent {\it Acknowledgements.}
This work has been carried out in the framework of the Labex Archimède (ANR-11-LABX-0033) and 
of the A*MIDEX project (ANR-11-IDEX-0001-02), funded by the ``Investissements d'Avenir" 
French Government programme managed by the French National Research Agency.

\bigskip 

{\scriptsize
{\sc Institut de Math\'ematiques de Marseille,  UMR 7373 CNRS/Aix-Marseille Universit\'e, 
Technop\^ole Ch\^ateau-Gombert, 
39 rue Fr\'ed\'eric Joliot-Curie,
13453 MARSEILLE Cedex 13,
France}

{\it Email address}:  {\tt laurent.manivel@math.cnrs.fr}
}
\end{document}